\documentclass[12pt]{amsart}

\usepackage{amssymb}
\usepackage{amscd}

\newcommand{\bt}{\begin{theorem}}
\newcommand{\et}{\end{theorem}}
\newcommand{\bp}{\begin{proposition}}
\newcommand{\ep}{\end{proposition}}
\newcommand{\bl}{\begin{lemma}}
\newcommand{\el}{\end{lemma}}
\newcommand{\br}{\begin{result}}
\newcommand{\er}{\end{result}}
\newcommand{\be}{\begin{equation}}
\newcommand{\ee}{\end{equation}}
\newcommand{\bc}{\begin{corollary}}
\newcommand{\ec}{\end{corollary}}
\newcommand{\bex}{\begin{example}}
\newcommand{\eex}{\end{example}}
    
\newtheorem{theorem}{Theorem}[section]
\newtheorem{corollary}[theorem]{Corollary}
\newtheorem{lemma}[theorem]{Lemma}
\newtheorem{proposition}[theorem]{Proposition}
\newtheorem{result}[theorem]{Result}
\newtheorem{example}[theorem]{Example}
    
\numberwithin{equation}{section}

\newcommand{\N}{\mathbb{N}}

\newcommand{\Z}{\mathbb{Z}}

\newcommand{\cB}{\mathcal{B}}
\newcommand{\cC}{\mathcal{C}}
\newcommand{\cD}{\mathcal{D}}
\newcommand{\cH}{\mathcal{H}}

\newcommand{\cJ}{\mathcal{J}}
\newcommand{\cK}{\mathcal{K}}

\newcommand{\cO}{\mathcal{O}}

\newcommand{\cPA}{\mathcal{P}\!\mathcal{A}}

\newcommand{\epr}{\hspace{\fill}$\Box$}

\newcommand{\bg}{\bigskip}
\newcommand{\med}{\medskip}
\newcommand{\sm}{\smallskip}

\font\caps=cmcsc10 scaled \magstep1   
    
\begin{document}
\noindent {\em Acta Sci. Math. (Szeged)} (2021) {\bf 87}, 367 -- 379\\
doi: 10.14232/actasm-020-558-7
\vspace{0.02in}\\
{\em arXiv version}: layout, fonts, pagination and numbering of lemmas, theorems and formulas may differ from those in the ACTASM published paper; moreover, references [12] and [13] in the published paper were printed with errors which have been corrected in this version (the author was not responsible for those errors -- according to one of the Associate Editors of Acta Sci. Math. (Szeged), they were caused by a mistake during their copy-editing process), and page numbers in reference [11] have also been corrected in this version.
\vspace{0.04in}\\
\title[On lattice isomorphisms of orthodox semigroups]{On lattice isomorphisms of orthodox semigroups}
\author[Simon M. Goberstein]{Simon M. Goberstein}\thanks{
Department of Mathematics and Statistics,
California State University, Chico, CA 95929-0525, U.S.A.\\
\indent \;\,e-mail: sgoberstein@csuchico.edu}
\begin{abstract}
Two semigroups are lattice isomorphic if the lattices of their subsemigroups are isomorphic, and a class of semigroups is lattice closed if it contains every semigroup which is lattice isomorphic to some semigroup from that class. An orthodox semigroup is a regular semigroup whose idempotents form a subsemigroup. We prove that the class of all orthodox semigroups in which every nonidempotent element has infinite order is lattice closed. 
\vspace{0.15in}\\
\noindent 2020 Mathematics Subject Classification: primary 20M15, 20M18, 20M19; secondary 08A30.
\end{abstract}
\maketitle
\font\caps=cmcsc10 scaled \magstep1   % Theorems, etc.
\def\bfseries{\normalsize\caps}
%\setcounter{section}{-1}
%\vspace{-0.2in}
\section{Introduction}
\med
Let $S$ be a semigroup. The set of all subsemigroups of $S$ (including, by convention, the empty one) is a lattice under set-theoretic inclusion, and the relationship between the properties of this lattice and the properties of $S$ has been studied in numerous publications for over 60 years. We say that two semigroups are lattice isomorphic whenever their subsemigroup lattices are isomorphic. If $S$ is isomorphic or antiisomorphic to any semigroup which is lattice isomorphic to $S$, then $S$ is called lattice determined in the class of all semigroups, and if for each isomorphism $\Phi$ of the subsemigroup lattice of $S$ onto that of a semigroup $T$, there exists an isomorphism or an antiisomorphism $\varphi$ of $S$ onto $T$ such that $U\Phi=U\varphi$ for every subsemigroup $U$ of $S$, then $S$ is termed strongly lattice determined. A class of semigroups is said to be lattice closed (in the class of all semigroups) if it contains every semigroup which is lattice isomorphic to some semigroup from that class. Of course, a lattice closed class of semigroups may contain members which are not lattice determined.  
\vspace{0.05in}\\
\indent
Finding lattice closed classes of semigroups and identifying lattice determined semigroups in such classes are among the most important problems in this area of research (see \cite{key20}). It is well known that the class of all groups is not lattice closed. However, the classes of all torsion-free groups and of all completely simple semigroups, which are neither groups nor left or right zero semigroups, are lattice closed \cite[\S\S\,34,\,38]{key20}, and it is also easily seen (and well known) that the class of bands is lattice closed. All semigroups mentioned in the preceding sentence are regular (in the sense of von Neumann for rings). A regular semigroup is orthodox if the set of its idempotents is a subsemigroup, and orthodox semigroups form one of the most important classes of regular semigroups.  Adopting for semigroups the standard group-theoretic term, we will say that a semigroup is torsion-free if all of its nonidempotent elements have infinite order. Note that according to Lemma \ref{205} in Section 2, a monogenic inverse semigroup whose generators have infinite order is torsion-free. The principal goal of this paper is to prove that the class of all torsion-free orthodox semigroups is lattice closed (see Theorem \ref{303} in Section 3). An important role in that proof is played by the author's theorem stating that every bisimple orthodox semigroup generated by a pair of mutually inverse elements of infinite order is strongly lattice determined \cite[Theorem 4.1]{key6}. Using Theorem \ref{303} (and its proof) together with Lemma \ref{205}, we also show that the class of all monogenic orthodox semigroups whose generators have infinite order is lattice closed and every semigroup in that class is torsion-free (Proposition \ref{308}). In Section 2, we establish a few auxiliary results and, for the reader's convenience, review some basic semigroup-theoretic concepts and include preliminaries on lattice isomorphisms of semigroups. 
\vspace{0.05in}\\
\indent
We use \cite{key1} and \cite{key11} as standard references for the algebraic theory of semigroups and, in general, follow the terminology and notation of these monographs. We also refer to \cite{key20} for information about lattice isomorphisms of semigroups and to \cite{key14} for a detailed description of the structure of monogenic inverse semigroups. 
\bg
\section{Preliminaries}
\med
We use the term ``order'' instead of ``partial order'' and refer to a linearly ordered set as a {\em chain}. In this paper, $\N$ stands for the set of all natural numbers (= positive integers), and $\N_0\,(=\N\cup\{0\})$ for the set of all nonnegative integers. Let $S$ be a semigroup. The equality (or ``diagonal'') relation on $S$ will be denoted by $\Delta_S$. We say that $x\in S$ is a {\em group element} of $S$ if it belongs to some subgroup of $S$; otherwise $x$ is a {\em nongroup element} of $S$. The set of nongroup elements of $S$ will be denoted by $N_S$, and the set of idempotents of $S$ by $E_S$. To indicate that $U$ is a subsemigroup of $S$, we write $U\leq S$. Under the convention that $\emptyset\leq S$, the set of all subsemigroups of $S$ ordered by inclusion is a lattice which we denote by $\operatorname{Sub}(S)$. If $w=w(x_1,\ldots,x_n)$ is a word in the alphabet $\{x_1,\ldots,x_n\}\subseteq S$, we will say shortly that $w$ is a word in $x_1,\ldots, x_n$, and if no confusion is likely, $w$ will be identified with its value in $S$.  As usual, $\langle X\rangle$ stands for the subsemigroup of $S$ generated by $X\subseteq S$, and if $X=\{x_1,\ldots,x_n\}$ is a finite subset of $S$, then $\langle X\rangle$ is written as $\langle x_1,\ldots, x_n\rangle$. Let $x\in S$. Then $\langle x\rangle$ is the {\em cyclic} subsemigroup of $S$ generated by $x$. If $\langle x\rangle$ is a finite semigroup, the number of its elements is called the {\em order} of $x$ and denoted by $o(x)$. If $\langle x\rangle$ is infinite, we write $o(x)=\infty$ and say that $x$ has {\em infinite order}. If every nonidempotent element of $S$ has infinite order, we call $S$ a {\em torsion-free} semigroup.  
\vspace{0.05in}\\
\indent If $S$ and $T$ are semigroups such that $\operatorname{Sub}(S)\cong \operatorname{Sub}(T)$, then $S$ and $T$ are called {\em lattice isomorphic}, and any isomorphism of $\operatorname{Sub}(S)$ onto $\operatorname{Sub}(T)$ is referred to as a {\em lattice isomorphism} of $S$ onto $T$. If $\Phi$ is a lattice isomorphism of $S$ onto $T$, then $\Phi$ is said to be {\em induced} by a mapping $\varphi \colon S\rightarrow T$ if $U\Phi=U\varphi$ for all $U\leq S$. If $S$ is isomorphic or antiisomorphic to any semigroup that is lattice isomorphic to it, then $S$ is called {\em lattice determined}, and if each lattice isomorphism of $S$ onto a semigroup $T$ is induced by an isomorphism or an antiisomorphism of $S$ onto $T$, we say that $S$ is {\em strongly lattice determined}. A class $\cK$ of semigroups is said to be {\em lattice closed} (in the class of all semigroups) if $\cK$ contains every semigroup which is lattice isomorphic to some semigroup that belongs to $\cK$. It is easily seen (and well known) that bands, in general, are not lattice determined -- for instance, any left or right zero semigroup $S$ is clearly lattice isomorphic to a chain (that is, a semilattice with a linear natural order) having the same cardinality as $S$. On the other hand, as an immediate corollary of \cite[Proposition 1]{key16} (reproduced as \cite[Proposition 3.2(a)]{key20}) we have: 
    
\br\label{201}  {\rm (From \cite[Proposition 1]{key16} or \cite[Proposition 3.2(a)]{key20})} A semigroup which is lattice isomorphic to a band is itself a band, that is, the class of all bands is lattice closed.
\er
    
For convenience of reference, we also record here the following two results (which imply, in particular, that the class of torsion-free semigroups is lattice closed). 

\br\label{202}{\rm \cite[Lemma 34.8]{key20}} The infinite cyclic group is strongly lattice determined.
\er
 
\br\label{203}{\rm (From \cite[Subsections 31.1--31.5]{key20})} Let $S$ and $T$ be torsion-free semigroups, and let $\Phi$ be a lattice isomorphism of $S$ onto $T$. Then $\Phi$ is induced by a unique bijection $\varphi \colon S\rightarrow T$ defined by the formula $\langle x\rangle\Phi=\langle x\varphi\rangle$ for all $x\in S$, and $(x^n)\varphi=(x\varphi)^n$ for all $x\in S$ and $n\in\N$. Thus the infinite cyclic semigroup is strongly lattice determined.
\er 
If $S$ and $T$ are torsion-free semigroups and $\Phi$ is a lattice isomorphism of $S$ onto $T$, the bijection $\varphi \colon S\rightarrow T$ in Result \ref{203} will be called the {\em $\Phi$-associated bijection} of $S$ onto $T$.
\vspace{0.05in}\\
\indent
Let $S$ be a semigroup. An element $x\in S$ is {\em regular} if there is $x'\in S$ such that $xx'x=x$, and if all elements of $S$ are regular, $S$ is a {\em regular semigroup}. A $\cD$-class $D$ of a regular semigroup $S$ is said to be {\em combinatorial} if all subgroups of $S$ contained in $D$ are trivial (or, equivalently, if $H_x=\{x\}$ for all $x\in D$), and $S$ itself is {\em combinatorial} if all $\cD$-classes of $S$ are combinatorial, so $S$ is combinatorial precisely when $\cH=\Delta_S$. If $x, x'\in S$ satisfy $xx'x=x$, then by setting $y=x'xx'$, we have $xyx=x$ and $yxy=y$, in which case, $x$ and $y$ are said to be mutually {\em inverse}. As in \cite{key5}, we indicate that $x$ and $y$ are mutually inverse elements of $S$ by writing $x\perp y$, and the phrase ``$x\perp y$ in $S$'' will always mean that $x,y\in S$ and $x\perp y$. We also denote by $V_S(x)$ the set of all inverses of $x\in S$, so $y\in V_S(x)$ if and only if $x\perp y$ in $S$. In general, a regular element of $S$ may have more than one inverse but if each $x\in S$ has a {\em unique} inverse (usually denoted by $x^{-1}$), then $S$ is an {\em inverse semigroup}. If $S$ is an inverse semigroup and $x\in S$, the inverse subsemigroup $\langle x, x^{-1}\rangle$ of $S$ is called {\em monogenic}, and if $S=\langle x, x^{-1}\rangle$ then $S$ is a {\em monogenic inverse semigroup}.
\vspace{0.05in}\\
\indent The {\em bicyclic semigroup} $\cB(a,b)$ is often defined as a semigroup with identity $1$ generated by the two-element set $\{a,b\}$ and given by one defining relation $ab=1$ \cite[\S 1.12]{key1}. Of course, $\cB(a,b)$ can also be defined without mentioning the identity as a semigroup given by the following presentation: $\cB(a,b)=\langle a, b\;|\;aba=a, bab=b, a^2b=a, ab^2=b\rangle$. It is well known that $\cB(a,b)$ is a combinatorial bisimple inverse semigroup with $b=a^{-1}$, each of its elements has a unique representation in the form $b^ma^n$, where $m$ and $n$ are nonnegative integers (and $a^0=b^0=ab$), and the semilattice of idempotents of $\cB(a,b)$ is a chain: $ab>ba>b^2a^2>\cdots$. Observe that if $S=\langle x, x^{-1}\rangle$ is a monogenic inverse semigroup such that $xx^{-1}>x^{-1}x$, then $S=\cB(x,x^{-1})$, and if $x^{-1}x>xx^{-1}$, then $S=\cB(x^{-1},x)$. Let $\cC=\{(m,n)\colon m, n\in\N_0\}$ with multiplication defined as follows: $(m,n)(p,q)=(m+p-r, n+q-r)$ where $r=\min\{n,p\}$. Then $\cC$ is a semigroup and the mapping $(m,n)\to b^ma^n$ is an isomorphism of $\cC$ onto $\cB(a,b)$ (see \cite[\S 1.12]{key1}). Therefore $\cC$ can be viewed as another copy of the bicyclic semigroup (in \cite{key14}, the bicyclic semigroup was, in fact, defined as the semigroup $\cC$ described above).
\vspace{0.05in}\\
\indent
A comprehensive description of the structure of monogenic inverse semigroups is given in \cite[Chapter IX]{key14}. We will recall only a few basic facts about them (see \cite{key14} for more details). Note, in particular, that if $S$ is a monogenic inverse semigroup, then $\cD=\cJ$, so the order relation on the set of $\cJ$-classes of $S$ is actually an order relation on the set of its $\cD$-classes. Let $S=\langle x, x^{-1}\rangle$ be a monogenic inverse semigroup. If $xx^{-1}=x^{-1}x$, then $S$ is a cyclic group. As noted in the preceding paragraph, if $xx^{-1}<x^{-1}x$ or $x^{-1}x<xx^{-1}$, then $S$ is a bicyclic semigroup. Assume that $S$ is neither a cyclic group nor a bicyclic semigroup (that is, $xx^{-1}$ and $x^{-1}x$ are incomparable with respect to the natural order on $E_S$). Then the set of $\cD$-classes of $S$ is a nontrivial chain with the largest element $D_x$. Moreover, one of the following holds: (i) $S$ is a free monogenic inverse semigroup, or (ii) $S$ has a smallest proper ideal $K$ (the {\em kernel} of $S$), which is either a bicyclic semigroup or a cyclic group. In case (i), all $\cD$-classes of $S$ are combinatorial and form an infinite chain: $D_x>D_{x^2}>\cdots>D_{x^m}>\cdots$, where the $\cD$-class $D_{x^m}$ consists of $(m+1)^2$ elements for each $m\in\N$. In case (ii), there is $n\in\N$ such that the $\cD$-classes of $S$ form a finite chain: $D_x>D_{x^2}>\cdots>D_{x^n}>D_{x^{n+1}}=K$, where for each $1\leq m\leq n$, the $\cD$-class $D_{x^m}$ is combinatorial and consists of $(m+1)^2$ elements. If $S$ is either a cyclic group or a bicyclic semigroup, it is convenient (and common) to assume that $S$ itself is the kernel of $S$, and with this assumption in mind we can state that if $S$ is an arbitrary monogenic inverse semigroup, then $S$ is not combinatorial if and only if it has a kernel which is a nontrivial cyclic group (and is the only nontrivial $\cH$-class of $S$).
\vspace{0.05in}\\
\indent
If $S=\langle x, x^{-1}\rangle$ is a free monogenic inverse semigroup, each $s\in S$ can be uniquely written in the form $s=x^{-p}x^qx^{-r}$ where $q>0$ and $0\leq p, r\leq q$, and $x^{-p}x^qx^{-r}\in E_S$ if and only if $p+r=q$; this copy of the free monogenic inverse semigroup was denoted in \cite{key14} by $C_3$. If the set $\{((m,n), (p,q))\in \cC\times \cC\colon m+p=n+q>0\}$ is equipped with multiplication of the direct square $\cC\times\cC$ of the bicyclic semigroup $\cC$, one obtains another copy of the free monogenic inverse semigroup, denoted in \cite{key14} by $C_2$. Clearly, $((m,n), (p,q))\in E_{C_2}$ if and only if $m=n$ and $p=q$. By setting $x=((1,0),(0,1))$ in $C_2$, it is easily shown that $C_2$ is isomorphic to $C_3=\langle x, x^{-1}\rangle$ (and a few other isomorphic copies of the free monogenic inverse semigroup were described in \cite{key14}). We are going to use several results about monogenic inverse semigroups expressed in \cite{key14} in terms of $C_2$ and its congruences, and for the reader's convenience we will recall some relevant definitions and notations from \cite{key14}.
\vspace{0.05in}\\
\indent
Let $S=\{((m,n), (p,q))\in \cC\times \cC\colon m+p=n+q>0\}$ with multiplication of the direct product $\cC\times\cC$. Let $\rho$ be a congruence on $S$. The least $n\in\N$ such that $((n,n),(0,0))$ and $((n+1,n+1),(0,0))$ are $\rho$-related is denoted by $l(\rho)$; if such $n$ does not exist, put $l(\rho)=\infty$. Similarly, $r(\rho)$ is the least natural number $n$ such that $((0,0),(n,n))$ and $((0,0),(n+1,n+1))$ are $\rho$-related, and if no such $n$ exists, then $r(\rho)=\infty$. By \cite[Lemma IX.2.6]{key14}, if $l(\rho)<\infty$ and $r(\rho)<\infty$, then $l(\rho)=r(\rho)$. Let $k\in\N$. As above, set $x=((1,0),(0,1))$.  If $l(\rho)=r(\rho)=k$ and $\rho|_{\langle x\rangle}$ has infinitely many classes, then $\rho$ is said to have {\em type} $(k,\omega)$. If $l(\rho)=\infty$ and $r(\rho)=k$, the {\em type} of $\rho$ is defined to be $(k,\infty^+)$, and if $l(\rho)=k$ and $r(\rho)=\infty$, the {\em type} of $\rho$ is $(k,\infty^-)$. The {\em weight} of an element $u=((m,n),(p,q))$ of $S$ is the number $w(u)=m+p$. Now for all $u=((m,n),(p,q)), u'=((m'\!,n'),(p'\!,q'))\in S$, define
\vspace{0.01in}\\
 \[\;u\rho_{(k,\,\omega)}u' \;\,\Longleftrightarrow u=u'\text{ or }[w(u), w(u')\geq k\text{ and }m-n=m'-n'],\]
\vspace{-0.1in} 
\[u\rho_{(k,\,\infty^{+})}u' \!\Longleftrightarrow u=u'\text{ or }[w(u), w(u')\geq k\text{ and }(m,n)=(m'\!,n')],\]
\vspace{-0.1in} 
\[\!\!\!\!\!u\rho_{(k,\,\infty^{-})}u' \!\Longleftrightarrow u=u'\text{ or }[w(u), w(u')\geq k\text{ and }(p,q)=(p'\!,q')].\]
\vspace{0.01in}\\
Then $\rho_{(k,\,\omega)}$ is the unique congruence on $S$ of type $(k,\,\omega)$ \cite[Proposition IX.2.11]{key14}, and $\rho_{(k,\,\infty^+)}$ is the unique congruence on $S$ of type $(k,\,\infty^+)$ \cite[Proposition IX.2.10]{key14}; of course, by symmetry, $\rho_{(k,\,\infty^-)}$ is the unique congruence on $S$ of type $(k,\,\infty^-)$.   
\vspace{0.05in}\\
\indent
Let $Q$ be a semigroup with zero $0$. Denote $Q^{\ast}=Q\setminus\{0\}$. Let $T$ be a semigroup disjoint from $Q$, and let $\varphi\colon Q^{\ast}\to T$ be a partial homomorphism. For all $a,b\in T\cup Q^{\ast}$ define $a\circ b$ as follows: $a\circ b=a(b\varphi)$ if $a\in T$ and $b\in Q^{\ast}$; $a\circ b=(a\varphi)b$ if $a\in Q^{\ast}$ and $b\in T$; $a\circ b=(a\varphi)(b\varphi)$ if $a,b\in Q^{\ast}$ and $ab=0$ in $Q$; $a\circ b=ab$ if either $a,b\in T$ or $a,b\in Q^{\ast}$ and $ab\ne 0$ in $Q$. Then $T\cup Q^{\ast}$ is a semigroup called a {\em retract ideal extension} (or simply an {\em ideal extension}) of $T$ by $Q$ whose operation $\circ$ is {\em determined by the partial homomorphism $\varphi$} (see \cite[Section I.9]{key14}).  
 \vspace{0.05in}\\
\indent
Consider again $S=\{((m,n), (p,q))\in \cC\times \cC\colon m+p=n+q>0\}$ with multiplication of the direct product $\cC\times\cC$, so $S$ is the free monogenic inverse semigroup $C_2$. For any $n\in\N$, denote by $I_n$ the set $\{u\in C_2\colon w(u)\geq n\}$. Then $I_n$ is an ideal of $C_2$, and we let $M_n=C_2/I_n$. In the notation of this paragraph, we can state the following result. 
 
\br\label{204}{\rm  (From \cite[Theorem IX.3.4 (i), (ii)]{key14})}
Let $k\in\N$. Then
\vspace{0.04in}\\
\indent
{\rm(i)} $C_2/\rho_{(k,\,\infty^+)}$ {\rm[respectively, }$C_2/\rho_{(k,\,\infty^-)}${\rm]} is an ideal extension of the bicyclic semigroup $\cC$ by $M_k$ with the operation in $\cC\cup M^{\ast}_k$ determined by the partial homomorphism $\theta: M^{\ast}_k\to\cC$ where $\theta: ((m,n), (p,q))\mapsto (m,n)\;{\rm[respectively, }\;\theta: ((m,n), (p,q))\mapsto (p,q){\rm]}${\rm;}
\vspace{0.04in}\\
\indent
{\rm(ii)} $C_2/\rho_{(k,\,\omega)}$ is an ideal extension of the cyclic group $\Z$ by $M_k$ with the operation in $\Z\cup M^{\ast}_k$ determined by the partial homomorphism $\eta: ((m,n), (p,q))\mapsto m-n$ of $M^{\ast}_k$ to $\Z$. 
\er
Of course, in Case (ii) of Result \ref{204}, instead of $\Z$ one can use a multiplicative infinite cyclic group $G=\langle g\rangle$, so $C_2/\rho_{(k,\,\omega)}$ will be an ideal extension of $G$ by $M_k$ with the operation in $G\cup M^{\ast}_k$ determined by the partial homomorphism $\eta:  ((m,n), (p,q))\mapsto g^{m-n}$ of $M^{\ast}_k$ to $G$. Take any nonidempotent element $((m,n), (p,q))$ of $M^{\ast}_k$. According to Result \ref{204}, $((m,n), (p,q))\theta\notin E_{\cC}$ in both situations considered in Case (i), and if in Case (ii) instead of $\Z$ we use a multiplicative infinite cyclic group $G=\langle g\rangle$ with the identity $1_G$, then $((m,n), (p,q))\eta\ne 1_G$.   

\bl\label{205} Let $S=\langle x, x^{-1}\rangle$ be a monogenic inverse semigroup such that $o(x)=\infty$. Then every nonidempotent element of $S$ has infinite order, that is, $S$ is torsion-free.
\el
{\bf Proof.} If $S$ is a free monogenic inverse semigroup, \cite[Lemma IX.3.8]{key14} shows that every nonidempotent element of $S$ has infinite order (this also follows from \cite[Exercise IX.2.14(iii)]{key14}, which asserts that  $\langle s, s^{-1}\rangle$ is a free monogenic inverse subsemigroup of $S$ for every $s\in S\setminus E_S$). Suppose that $S$ is not free. Then $S$ has a kernel $K$ (perhaps coinciding with $S$), which is either a bicyclic semigroup or a cyclic group (the latter must be infinite since $o(x)=\infty$).  
\vspace{0.05in}\\
\indent
Take any $s\in S\setminus E_S$. If $s\in K$, it is immediate that $o(s)=\infty$. Now assume that $K$ is a proper ideal of $S$ and $s\notin K$. Thus there is $k\geq 2$ such that $D_x>\cdots>D_{x^{k-1}}>D_{x^k}=K$, and $s\in D_x\cup\cdots\cup D_{x^{k-1}}$. In the notation of the paragraph preceding Result \ref{204}, $K=I_k$ and $D_x\cup\cdots\cup D_{x^{k-1}}=M^{\ast}_k$. Moreover, if $K$ is an infinite cyclic group, then $S=C_2/\rho_{(k,\,\omega)}$, and if $K$ is a bicyclic semigroup, then $S=C_2/\rho_{(k,\,\infty^+)}$ or $S=C_2/\rho_{(k,\,\infty^-)}$. Since all nonidempotent elements of a free monogenic inverse semigroup have infinite order, $\langle s\rangle\nsubseteq M^{\ast}_k$ and hence $\langle s\rangle\cap K\ne\emptyset$.  Let $l$ be the smallest natural number such that $s^l\in K$. Then $s^{l-1}\in M^{\ast}_k$. Using Result \ref{204}, the remark that follows it, and the formulas for the powers of a nonidempotent element $((m,n), (p,q))$ of $C_2$ given in \cite[Lemma IX.3.8]{key14}, one can deduce that if $K$ is a bicyclic semigroup, then $s^l=(s\theta)(s^{l-1}\theta)\notin E_S$, and if $K$ is an infinite cyclic group, then $s^l=(s\eta)(s^{l-1}\eta)\notin E_S$. We have shown that in both cases $s^l\in K\setminus E_S$. If $o(s)<\infty$, the cyclic semigroup $\langle s\rangle$ is finite and since $s^l\in\langle s\rangle$, we would have $o(s^l)<\infty$, contradicting the fact that all nonidempotent elements of $K$ have infinite order. Therefore $o(s)=\infty$. The proof is complete.
 \epr 
\vspace{0.05in}\\
\indent
 An {\em orthodox semigroup} is a regular semigroup in which the idempotents form a subsemigroup. By \cite[Theorem VI.1.1]{key11}, if $S$ is a regular semigroup, the following conditions are equivalent: (a) $S$ is orthodox; (b) $V_S(e)\subseteq E_S$ for all $e\in E_S$; (c) $V_S(b)V_S(a)\subseteq V_S(ab)$ for all $a,b\in S$. It follows that if $a\perp b$ in an orthodox semigroup $S$, then $a\in N_S$ if and only if $b\in N_S$, and $a^n\perp b^n$ for all $n\in\N$, which implies that $o(a)=o(b)$. According to the terminology introduced by the author in \cite{key5} (by analogy with the inverse semigroup case), a {\em monogenic orthodox semigroup} is an orthodox semigroup generated by a pair of mutually inverse elements. In what follows, the phrase ``let $A=\langle a,b\rangle$ be a monogenic orthodox semigroup'' will always mean that $A$ is an orthodox semigroup with $a\perp b$ in $A$. 
\vspace{0.05in}\\
\indent
It was shown by Hall (see \cite[Theorem 2]{key9} or \cite[Theorem VI.1.10]{key11}) that a regular semigroup $S$ is orthodox if and only if for all $x, y\in S$, if $V_S(x)\cap V_S(y)\ne\emptyset$ then $V_S(x)= V_S(y)$. Let $S$ be an arbitrary orthodox semigroup. Since $E_S$ is a band, it is a semilattice $Y$ of rectangular bands $E_{\alpha}\,({\alpha}\in Y)$, that is, $E_S=\bigcup_{\alpha\in Y}E_{\alpha}$ where $Y$ is a semilattice and $E_{\alpha}\,(\alpha\in Y)$ are pairwise disjoint rectangular bands satisfying $E_{\alpha}E_{\beta}\subseteq E_{\alpha\beta}$ for all $\alpha,\beta\in Y$ (see \cite[\S 3]{key9} or \cite[Section VI.1]{key11}); as in \cite{key9} and \cite{key11}, we may denote $E_{\alpha}$ by $E(e)$ if $e\in E_{\alpha}$. Now let $\gamma_S=\{(x,y)\in S\times S\colon V_S(x)=V_S(y)\}$. If no confusion is likely, we will omit the subscript $S$ in $\gamma_S$ and in $V_S(x)$ for all $x\in S$. It was proved by Hall that $\gamma$ is the smallest inverse semigroup congruence on $S$ (see \cite[Theorem 3]{key9} or \cite[Theorem VI.1.12]{key11}), so $S/\gamma\,(=S\gamma^{\natural})$ is the maximum inverse semigroup homomorphic image of $S$. According to \cite[Remark 1]{key9}, $e\gamma=V(e)=E(e)$ for all $e\in E_S$. Thus $e\gamma$ is a rectangular band for every $e\in E_S$, the semilattices $E_{S/\gamma}$ and $Y$ are isomorphic, and $s\gamma^{\natural}\in E_{S/\gamma}$ if and only if $s\in E_S$ for all $s\in S$. It follows that for any $x\in S$, we have $o(x)=\infty$ if and only if $o(x\gamma)=\infty$. Therefore $S$ is torsion-free if and only if $S/\gamma$ is torsion-free. {\em The results reviewed in this paragraph will be used below without any further reference or explanation.}
\vspace{0.05in}\\
\indent
Let $S$ be an orthodox semigroup. Since $\gamma\cap\cH=\Delta_S$ (see \cite[Chapter VI, formula (1.16)]{key11}), if $H$ is any $\cH$-class of $S$, the restriction of $\gamma^{\natural}$ to $H$ is an injection. As noted in \cite{key3}, from this observation and the fact that a monogenic inverse semigroup contains a nontrivial $\cH$-class if and only if it has a kernel which is a nontrivial cyclic group (and is the only nontrivial $\cH$-class of that semigroup), one can deduce the following assertion:
\setcounter{theorem}{5} 
 \br\label{206}{\rm (From \cite[Proof of Theorem 2.1]{key3}} Let $A=\langle a, b\rangle$ be a monogenic orthodox semigroup. If $H$ is a nontrivial $\cH$-class of $A$, then $A$ has a kernel $K$ and $H\subseteq K$.
 \er

Let us show that the assertion stated in Result \ref{206} can be made more precise.

\bp\label {207}
Let $A=\langle a, b\rangle$ be a monogenic orthodox semigroup. Suppose that $H$ is a nontrivial $\cH$-class of $A$. Then $H$ is a cyclic group isomorphic to the kernel of $A/\gamma$ and $A$ has a completely simple kernel $K$, which contains $H$ and is isomorphic to the direct product of $H$ and a rectangular band.
\ep
{\bf Proof.} As noted prior to Result \ref{206}, the restriction of $\gamma^{\natural}$ to $H$ is an injection from $H$ to the kernel of $A/\gamma$, which is a nontrivial cyclic group. Denote the kernel of $A/\gamma$ by $G$, and let $K=G\left(\gamma^{\natural}\right)^{-1}$. Then $K$ is an ideal of $A$ and $H\subseteq K$. Clearly, $E_K=\{e\in K\colon e\gamma^{\natural}=1_G\}$ where $1_G$ is the identity of $G$. If $e, f\in E_K$, then $e\gamma^{\natural}=f\gamma^{\natural}$, so $e\perp f$ and hence $(e, f)\in\cD$. Therefore $K$ is a regular $\cD$-class (= $\cJ$-class) of $A$, which shows that $K$ is the kernel of $A$. If $e$ is any element of $E_K$, then $E_K=\{f\in K\colon f\in e\gamma\}=e\gamma$. Thus $E_K$ is a rectangular band. By \cite[Corollary IV.3.5]{key13} (or \cite[Exercise III.12]{key11}), $K\cong H\times E_K$. Since each $e\in E_K$ is primitive, $K$ is completely simple. Clearly, the maximum inverse semigroup homomorphic image of $H\times E_K$ is isomorphic to $H$, and $K\gamma^{\natural}=G$, so $H\cong G$. The proof is complete.	\epr  
\med  
\section{The class of torsion-free orthodox semigroups is lattice closed}
\sm
It is well known that cyclic groups and the bicyclic semigroup are the only bisimple monogenic inverse semigroups. However, as shown in \cite{key5}, the class of bisimple monogenic orthodox semigroups is substantially more diverse. In particular, the author constructed in \cite{key5} a family of pairwise nonisomorphic bisimple orthodox semigroups $\cO_{(\nu,\,\mu)}(a,b)$ indexed by ordered pairs $(\nu,\mu)\in\N^{\ast}\times\N^{\ast}$ (where $\N^{\ast}=\N\cup\{\infty\}$), each being generated by a pair of mutually inverse elements $a$ and $b$ satisfying $ab=a^2b^2$ and $ba\neq b^2a^2$, and proved that if $S$ is an arbitrary bisimple monogenic orthodox semigroup with nongroup generators, then $S$ or its dual is isomorphic to one of the semigroups of that two-parameter family. 
    
\bt\label{301}{\rm (From \cite[Lemma 2.7 and Theorem 2.9]{key5})} Let $S$ be an orthodox semigroup, let $a$ be an arbitrary nongroup element of $S$, and let $b\in V(a)$. Then either $\{a, b, ab, ba\}$ is a $\cD$-class of $\langle a, b\rangle$ such that $\langle a,b\rangle\setminus\{a, b, ab, ba\}$ is an ideal of $\langle a, b\rangle$, or $\langle a, b\rangle$ is a bisimple monogenic orthodox semigroup and $o(a)=o(b)=\infty$, in which case $\langle a, b\rangle$ or its dual is isomorphic to $\cO_{(\nu,\,\mu)}(a,b)$ for some $\mu, \nu\in\N^{\ast}$.
\et
    
As shown in \cite{key18} (or \cite[Theorem 41.8]{key20}), the bicyclic semigroup is strongly lattice determined. A much more general result was established in \cite{key6} where it was proved that {\em for all $\mu, \nu\in\N^{\ast}$, the semigroups $\cO_{(\nu,\,\mu)}(a,b)$ are strongly lattice determined} (the bicyclic semigroup is just one member of that infinite two-parameter family -- namely, $\cO_{(1, 1)}(a,b)$ is bicyclic); this is the main part of the following theorem, from which it is obtained when the generators $a$ and $b$ of a bisimple monogenic orthodox semigroup $S=\langle a, b\rangle$ are assumed to be nongroup.
    
\bt\label{302} \cite[Theorem 4.1]{key6} Let $S=\langle a,b\rangle$ be an arbitrary bisimple monogenic orthodox semigroup such that $a$ (and hence $b$) has infinite order. Then $S$ is strongly lattice determined.
\et 
We are ready to establish the main result of the paper:     
\bt\label{303} Let $S$ be a torsion-free orthodox semigroup, and let $\Phi$ be a lattice isomorphism of $S$ onto a semigroup $T$. Then $T$ is also a torsion-free orthodox semigroup. Thus the class of all torsion-free orthodox semigroups is lattice closed.
\et
{\bf Proof.} By Result \ref{201}, $E_T\ne\emptyset$ and $E_T$ is a subsemigroup of $T$ such that $E_T=E_S\Phi$. Since the class of torsion-free semigroups is lattice closed, we can apply Result \ref{203} and denote by $\varphi$ the $\Phi$-associated bijection of $S$ onto $T$. Group elements of $T$ are regular. To prove that $T$ is an orthodox semigroup, it remains to show that all nongroup elements of $T$ are regular.
\vspace{0.05in}\\
\indent
Take an arbitrary $x\in N_T$, and let $a=x\varphi^{-1}$. By Result \ref{202}, the infinite cyclic group is strictly lattice determined. Since $S$ is torsion-free, if $a$ were a group element of $S$, it would be contained in some infinite cyclic subgroup $G$ of $S$, so $x$ would belong to the infinite cyclic subgroup $G\varphi$ of $T$, contradicting the assumption that $x\in N_T$. Therefore $a\in N_S$. Choose an arbitrary $b\in V(a)$. Then $\langle a, b\rangle$ is a monogenic orthodox semigroup such that $a,b\in N_S$ and $o(a)=o(b)=\infty$. Let $y=b\varphi$. By Result \ref{203}, $o(x)=o(y)=\infty$. Note that 
\[\langle a, b\rangle\Phi=\left(\langle a\rangle\vee\langle b\rangle\right)\Phi=\langle a\rangle\Phi\vee\langle b\rangle\Phi=\langle x\rangle\vee\langle y\rangle=\langle x, y\rangle. \]
\indent
Suppose that  $\langle a, b\rangle$ is a bisimple monogenic orthodox semigroup. Then, by Theorem \ref{302}, it is strongly lattice determined. Therefore $\langle x, y\rangle$ is also a bisimple monogenic orthodox semigroup which is isomorphic or antiisomorphic to $\langle a, b\rangle$. In particular, $x\perp y$ in $T$ and hence $x$ is a regular element of $T$.
\vspace{0.05in}\\
\indent 
Now assume that $\langle a, b\rangle$ is not bisimple. Then, according to Theorem \ref{301}, $\{a, b, ab, ba\}$ is a $\cD$-class of $\langle a, b\rangle$ such that $\langle a,b\rangle\setminus\{a, b, ab, ba\}$ is an ideal of $\langle a, b\rangle$. Let $I=\langle a,b\rangle\setminus\{a, b, ab, ba\}$ and $J=I\Phi$. Take an arbitrary $s\in\langle a, b\rangle$. Since $\langle s, I\rangle=\langle s\rangle\vee I$ and $\langle s\varphi, J\rangle=\langle s\varphi\rangle\vee J$, we conclude that  
\[ \langle s, I\rangle\Phi=\left(\langle s\rangle \vee I\right)\Phi=\langle s\rangle\Phi \vee I\Phi=\langle s\varphi \rangle \vee J=\langle s\varphi, J\rangle. \] 
If $s\in I$ then $\langle s, I\rangle = I=\{s\}\cup I$. If $s\in\{ab, ba\}$, it is also clear that $\langle s, I\rangle= \{s\}\cup I$. Finally, if $s\in\{a, b\}$ then again $\langle s, I\rangle= \{s\}\cup I$ because $s^k\in I$ for all $k\geq 2$. We have shown that $\langle s, I\rangle=\{s\}\cup I$, and hence $\langle s\varphi, J\rangle=\langle s, I\rangle\varphi=\left(\{s\}\cup I\right)\varphi=\{s\}\varphi\cup I\varphi=\{s\varphi\}\cup J$. Since $a, b\notin I$, we have $x, y\notin J$. Let $e=(ab)\varphi$ and $f=(ba)\varphi$. Then $e, f\in E_T$ and $e, f\notin J$ because $ab, ba\notin I$. Note that $\langle x, y\rangle=\langle a, b\rangle\Phi=\left(\{a, b, ab, ba\}\cup I\right)\varphi=\{a\varphi, b\varphi, (ab)\varphi, (ba)\varphi\}\cup I\varphi$. Therefore we have    
\setcounter{equation}{3}
\be\langle x, y\rangle=\{x, y, e, f\}\cup J\text{ and }\{x, y, e, f\}\cap J=\emptyset.
\label{eq:1st}
\ee
Since $\langle a\rangle\cap I=\langle a\rangle\setminus\{a\}$, it follows that 
\[ \langle x\rangle\cap J= \langle a\rangle\Phi\cap I\Phi=\left(\langle a\rangle\cap I\right)\Phi=\left(\langle a\rangle\setminus\{a\}\right)\Phi=\langle a\rangle\varphi\setminus\{a\}\varphi=\langle x\rangle\setminus\{x\}, \]
which shows that $x^k\in J$ for all $k\geq 2$. By symmetry, we also have $y^k\in J$ for all $k\geq 2$.
\setcounter{theorem}{4}
\bl\label{305} One of the following two statements is true:
\vspace{0.03in}\\
\indent
{\rm (i)} $e=xy$ and $f=yx$, or
\vspace{0.03in}\\
\indent
 {\rm (ii)} $e=yx$ and $f=xy$.
\el
{\bf Proof.} Since $e\in\langle x, y\rangle$, we can choose a shortest possible word $u$ in $x, y$ representing $e$. Thus $e=u$ and no word in $x, y$ which is shorter than $u$ has value $e$. 
\vspace{0.05in}\\
\indent
{\bf Case 1:} The first letter of $u$ is $x$.
\vspace{0.04in}\\
\indent Since $ab\ne a$, it follows that $e\ne x$. Therefore $u=xv$ for some nonempty word $v$ in $x, y$, and $e\ne v$ because $v$ is shorter than $u$. Since $e=xv$, we conclude that $e\in\langle x, v\rangle=\langle x\rangle\vee\langle v\rangle$. Let $c=v\varphi^{-1}$. Then
\[ab=e\varphi^{-1}\in\left(\langle x\rangle\vee\langle v\rangle\right)\Phi^{-1}=\langle x\rangle\Phi^{-1}\vee\langle v\rangle\Phi^{-1}=\langle x\varphi^{-1}\rangle\vee\langle v\varphi^{-1}\rangle=\langle a\rangle\vee\langle c\rangle=\langle a, c\rangle.\]
If $c\in I$, then $ab\in\langle a, I\rangle=\{a\}\cup I$, which is not true. Hence $c\in\{a, b, ab, ba\}$. If $c=a$, then $v=a\varphi=x$, so that $e=x^2\in J$; a contradiction. Since $v\ne e$, it is clear that $c\ne ab$. If $c=ba$, then $ab\in\langle a, ba\rangle$, contradicting the easily verified fact that $\langle a, ba\rangle\setminus I=\{a, ba\}$. Since $c\notin\{a, ab, ba\}$, it follows that $c=b$, so $v=b\varphi=y$. Therefore $e=xy$. 
\vspace{0.05in}\\
\indent
{\bf Case 2:} The first letter of $u$ is $y$.
\vspace{0.04in}\\
\indent By a symmetric argument to that used in Case 1, we deduce that $e=yx$.
\vspace{0.05in}\\
\indent Similarly to the above, starting with $f\in\langle x, y\rangle$, we prove that either $f=xy$, or $f=yx$. Since $e\ne f$, note that if $e=xy$ then $f\ne xy$ and hence $f=yx$, whereas if $e=yx$ then $f\ne yx$ and so $f=xy$.  Therefore either (i) $e=xy$ and $f=yx$, or (ii) $e=yx$ and $f=xy$. This completes the proof of Lemma \ref{305}. \epr
\vspace{0.05in}\\
\indent 
By Lemma \ref{305}, $xy$ and $yx$ are distinct idempotents of $T$ such that $\{xy, yx\}=\{e, f\}$. Therefore, using (\ref{eq:1st}), we conclude that
\setcounter{equation}{5}
\be\langle x, y\rangle=\{x, y, xy, yx\}\cup J\text{ and }\{x, y, xy, yx\}\cap J=\emptyset.
\label{eq:2nd}
\ee
\noindent
Our goal is to prove that $x$ is a regular element. In fact, we are going to show that $x\perp y$. 
\vspace{0.05in}\\
\indent
Since $xyx\in\langle x, y\rangle$, according to (\ref{eq:2nd}), we have
\be xyx\in \{x, y, xy, yx\}\cup J.
\label{eq:3rd}
\ee
\noindent 
If $xyx\in J$, then $xy=(xy)^2=(xyx)y\in\langle y, J\rangle=\{y\}\cup J$; a contradiction because $xy\ne y$ and $xy\notin J$. Therefore $xyx\notin J$. If $xyx=y$, then $xy=(xyx)y=y^2\in J$, which is not true. Thus $xyx\ne y$. Assume that $xyx=xy$. Then $xy=(xyx)y=xy^2\in\langle x, J\rangle=\{x\}\cup J$; a contradiction since $xy\ne x$ and $xy\notin J$. Hence $xyx\ne xy$. By a dual argument, $xyx\ne yx$. We have shown that $xyx\notin\{y, xy, yx\}\cup J$. In view of (\ref{eq:3rd}), it follows that $xyx=x$, that is, $x$ is a regular element of $T$. By symmetry, we also have $yxy=y$, so $x\perp y$.  The proof of Theorem \ref{303} is complete.  \epr
\vspace{0.07in}\\
{\bf Remark.} By \cite[Theorem]{key12} and comments after \cite[Result 4.5]{key4}, if $A\!=\!\langle a, a^{-1}\rangle$ is a monogenic inverse semigroup, which is neither bicyclic nor a group, and $\Psi$ is an isomorphism of the partial automorphism monoid of $A$ onto that of a semigroup $B$, then $B\!\cong\! A$. One step in the proof of this is \cite[Lemma 1]{key12} according to which $B\!=\!\langle x, y\rangle$ for mutually inverse $x$ and $y$ defined by $\Delta_{\langle x\rangle}\!=\!\Delta_{\langle a\rangle}\Psi$ and $\Delta_{\langle y\rangle}\!=\!\Delta_{\langle a^{-1}\rangle}\Psi$. The assumption that $A$ is a monogenic {\em inverse} semigroup (neither bicyclic nor a group) and {\em partial automorphism monoids} of $A$ and $B$ are isomorphic, is crucial for all results of \cite{key12}. However, due to Theorem \ref{301}, certain arguments in the proof of \cite[Lemma 1]{key12} have natural analogues in the more general setting of this paper. For instance, the argument in the last paragraph of the proof of Theorem \ref{303} is similar to the one in the final paragraph of the proof of \cite[Lemma 1]{key12}.  
\vspace{0.07in}\\
\indent 
Using Lemma \ref{205} and Theorem \ref{303} (and its proof), we can establish the following result.
\setcounter{theorem}{7}
\bp\label{308} Let $A=\langle a, b\rangle$ be a monogenic orthodox semigroup such that $a$ (and hence $b$) has infinite order. Then $A$ is torsion-free and any semigroup lattice isomorphic to $A$ is a torsion-free monogenic orthodox semigroup. Thus the class of all monogenic orthodox semigroups whose generators have infinite order is lattice closed and every semigroup in that class is torsion-free.
\ep
{\bf Proof.} Since $A/\gamma=\langle a\gamma, b\gamma\rangle$ is a monogenic inverse semigroup with $b\gamma=(a\gamma)^{-1}$ such that $o(a\gamma)=\infty$, according to Lemma \ref{205}, $A/\gamma$ is torsion-free. Therefore  $A$ is torsion-free.
\vspace{0.05in}\\
\indent
By Result \ref{203}, $\langle a\rangle$ and $\langle b\rangle$ are strongly lattice determined. Thus $\langle a\rangle\Phi=\langle x\rangle$ and $\langle b\rangle\Phi=\langle y\rangle$ for some $x, y\in X$ such that $o(x)=o(y)=\infty$. Then 
\[X=A\Phi=\langle a, b\rangle\Phi=\left(\langle a\rangle\vee\langle b\rangle\right)\Phi=\langle a\rangle\Phi\vee\langle b\rangle\Phi=\langle x\rangle\vee\langle y\rangle=\langle x, y\rangle.\]
If $A$ is bisimple, according to Theorem \ref{302}, $\langle x, y\rangle$ is a bisimple monogenic orthodox semigroup which is isomorphic or antiisomorphic to $\langle a, b\rangle$. Suppose that $A$ is not bisimple. Then, as shown in the proof of Theorem \ref{303}, $x\perp y$. Moreover, since $A$ is torsion-free, applying Theorem \ref{303}, we conclude that $X$ is a torsion-free orthodox semigroup. Therefore in all cases $X=\langle x, y\rangle$ is a torsion-free monogenic orthodox semigroup. This completes the proof.  \epr
\vspace{0.5in}\\ 
\begin{center}
{\bf Acknowledgement}
\end{center}
\sm The author would like to thank a careful referee for useful comments, which have helped shorten the proofs of Propositions \ref{207} and \ref{308}, resulting in an improved version of the paper.  
\vspace{0.05in}\\


\begin{thebibliography}{99}
\bg
\bibitem{key1} A.~H.~Clifford and G.~B.~Preston, {\it The Algebraic Theory of Semigroups}, Math. Surveys No. 7, Amer. Math. Soc., Providence, R.I.; Vol. I, 1961; Vol. II, 1967.
\bibitem{key3} C.~Eberhart and W.~Williams, Elementary orthodox semigroups, {\it Semigroup Forum}, $\mathbf {29}$ (1984), 351--364.
\bibitem{key4} S.~M.~Goberstein, $\cPA$-isomorphisms of inverse semigroups, {\it Algebra Universalis}, $\mathbf {53}$ (2005), 407--432.
\bibitem{key5} S.~M.~Goberstein, Bisimple monogenic orthodox semigroups, {\it Semigroup Forum}, $\mathbf {78}$ (2009), 310--325.
\bibitem{key6} S.~M.~Goberstein, Lattice isomorphisms of bisimple monogenic orthodox semigroups, {\it Semigroup Forum}, $\mathbf {83}$ (2011), 250--280.
\bibitem{key9} T.~E.~Hall, On regular semigroups whose idempotents form a subsemigroup, {\it Bull. Australian Math. Soc.}, $\mathbf {1}$ (1969), 195--208. 
\bibitem{key11} J.~M.~Howie, {\it An Introduction to Semigroup Theory}, Academic Press, London, 1976.
\bibitem{key12}  A.~L.~Libih, Local automorphisms of monogenic inverse semigroups,
{\it Theory of Semigroups and Its Applications}, No. 4, Saratov Univ. Press, 1978, 54--59 (in Russian).
\bibitem{key13} M.~Petrich, {\it Introduction to Semigroups}, Merrill, Columbus, Ohio, 1973.
\bibitem{key14} M.~Petrich, {\it Inverse Semigroups}, Wiley, New York, 1984.
\bibitem{key16}  L.~N.~Shevrin, Lattice properties of idempotent semigroups. I, {\it Sib. Mat. J.}, $\mathbf {6}$ (1965), 459--474 (in Russian).
\bibitem{key18}  L.~N.~Shevrin, The bicyclic semigroup is determined by its subsemigroup lattice, {\it Simon Stevin}, $\mathbf {67}$ (1993) Supplement (December 1993), 49--53.
\bibitem{key20}  L.~N.~Shevrin and A.~J.~Ovsyannikov, {\it Semigroups and Their Subsemigroup Lattices}, Kluwer Academic Publishers, Dordrecht, 1996.
\end{thebibliography}
\end{document}